\documentclass [a4paper,twoside,12pt]{article}
\usepackage[latin1]{inputenc}
\usepackage{amsmath,amsfonts,amssymb}
\usepackage{vmargin,graphicx,theorem}
\usepackage[english]{babel}
\usepackage{enumerate}
\usepackage{color}
\usepackage[colorlinks=true]{hyperref}

\setpapersize[portrait]{A4}
\setmarginsrb{2cm}{2.5cm}{2cm}{2cm}{0cm}{0.1cm}{0.5cm}{2cm}

\newcommand{\red}{\color{red} }

\newcommand{\dP}{\ensuremath{\mathbb{P}}}

\newcommand{\dR}{\ensuremath{\mathbb{R}}}

\newtheorem{ethm}{Theorem}[section]
\newtheorem{fthm}{Th\'{e}or\`{e}me}[section]
\newtheorem{ecor}[ethm]{Corollary}

\newtheorem{eprop}[ethm]{Proposition}

\newtheorem{elem}[ethm]{Lemma}
\newtheorem{flem}[fthm]{Lemme}

\newtheorem{erem}[ethm]{Remark}

\newcommand{\proofend}{~$\rhd$}
\newcommand{\proofbegin}{~$\lhd$}

\renewcommand{\phi}{{q}}

\renewcommand{\geq}{\geqslant}

\newcommand{\N}{\ensuremath{\mathbb{N}}}
\newcommand{\R}{\dR}

\begin{document}

\title{\sl Super-Poincar\'e and Nash-type inequalities  for Subordinated Semigroups}
\author{Ivan Gentil
\thanks{Institut Camille Jordan,
Universit\'e Claude Bernard Lyon 1.
43 boulevard du 11 novembre 1918.
69622 Villeurbanne cedex.
France. {\em gentil@math.univ-lyon1.fr }},
 Patrick Maheux
 \thanks{F\'ed\'eration Denis Poisson, D\'epartement de Math\'ematiques,
Universit\'e  d'Orl\'eans.
Rue de Chartres.
B.P. 6759 - 45067 Orl\'eans cedex 2.
France. \em{patrick.maheux@univ-orleans.fr}}
}

\date{\today}
\maketitle
\thispagestyle{empty}
\abstract{
We prove that if  a super-Poincar\'e inequality is  satisfied by     an infinitesimal   generator  $-A$ of a symmetric contracting semigroup 
then it   implies a corresponding     super-Poincar\'e inequality for $-g(A)$ with  any Bernstein function $g$. We also study the converse statement.
We deduce similar results for the  Nash-type inequality. 
Our results   applied  to fractional powers of $A$ and to  $\log(I+A)$   and thus generalize some  results of   \cite{bm} and
\cite{w1}.   We provide several  examples.}
\\

\noindent
{\bf Key words:} 
Super-Poincar\'e inequality, Nash-type inequality, Symmetric semigroup, Subordination in the sense of Bochner, Bernstein function, Super-Poincar\'e profile.
\\

\noindent
{\bf AMS  2010}: 39B62, 47D60, 26A33, 46T12. 
 

\section{Introduction and main results}\label{intro}

 Let $(T_t)_{t>0}$ be a strongly continuous semigroup  on $L^2(X,\mu)$ with 
$(X,\mu)$  a $\sigma$-finite measure space. We assume that, for each $t>0$, $T_t$  is a symmetric   contraction on $L^2$.
The infinitesimal generator of $(T_t)_{t>0}$  on $L^2$ denoted by $-A$ is a non-negative, closed and symmetric operator densely defined  on $L^2$. 
We shall not distinguish $-A$ and    its  self-adjoint Friedrich extension.
Moreover, we suppose  that, for each $t>0$,  $T_t$  is a contraction from  $L^1\cap L^2$ to  itself. So, for each $t>0$, the operator $T_t$ can be uniquely extended as a contraction on each  $L^p$, $1\leq p\leq +\infty$. This extension $T_t^{(p)} $ will be  again  denoted by $T_t$.
Recall that the infinitesimal generator $-A$  on $L^2(\mu)$ of  $(T_t)_{t>0}$  is defined by
$$
-Af=\lim_{t\rightarrow 0^+} \frac{T_tf-f}{t}\in L^2
$$
on the domain ${\mathcal D}(A)$ which is the set of functions $f$ such that the limit just above exists in $L^2(\mu)$.
We denote by $(., .)$ the inner product on $L^2$ and by $\vert\vert. \vert\vert_p$ the $L^p$-norm.


We shall say that $A$ satisfies {\it a  super-Poincar\'e inequality}   with rate function $\beta:(0,+\infty) \rightarrow  (0,+\infty)$ if,  for any $f\in {\mathcal D}(A)\cap L^1(\mu)$,
 \begin{equation}\label{poini}
\vert\vert f\vert \vert_2^2\leq r(Af,f)+\beta(r)\vert\vert f\vert \vert_1^2, \; r>0.
 \end{equation}
More generally, we say that $A$ satisfies a $(r_0,r_1)$-super-Poincar\'e inequality if 
(\ref{poini}) holds for $r\in(r_0,r_1)$ with  $0< r_0<r_1\leq +\infty$.
See examples below  Theorem \ref{main2}.
Note  that we can always assume that $\beta$ is non-increasing by considering what we shall call  the 
{\em super-Poincar\'e profile} $\beta_p$ ($\leq \beta$)
of $A$ defined, for any $r>0$,  by
\begin{equation}\label{profile}
\beta_p(r):=\sup\{ \vert\vert f \vert\vert_2^2-r(Af,f): f\in {\mathcal D}(A)\cap L^1(\mu), \vert\vert f\vert\vert_1\leq 1\}.
\end{equation}

We shall say that $A$ satisfies a {\it Nash-type inequality} with  non-decreasing rate function
 $D:[0,+\infty)\rightarrow [0,+\infty)$ if, for any $f\in {\mathcal D}(A)\cap L^1(\mu)$,
\begin{equation}\label{nashi}
\vert\vert f\vert \vert_2^2\, D(\vert\vert f\vert \vert_2^2)\leq (Af,f), \quad \vert\vert f\vert \vert_1\leq 1.
\end{equation}
 $D$  is then a {\it Nash function}  for $A$.
It is well known that the inequalities (\ref{poini})  and  (\ref{nashi}) are essentially equivalent,
see Proposition \ref{nashpoin}.

We now briefly recall some   definitions and some facts   about  the subordination of semigroups in the sense of Bochner.

A Bernstein function  $g$   is a    $\mathcal C^\infty$ function $g:(0,+\infty)\rightarrow (0,+\infty)$ satisfying 
$$
(-1)^{n-1}g^{(n)}(s)\geq 0
$$ for any $n\in \N^*$, $s>0$. 
There exists a convolution semigroup of sub-probability measures $(\nu^g_t)_{t>0}$ on $(0,+\infty)$ with density
$(\eta^g_t)_{t>0}$ with respect to the Lebesgue measure $ds$
 such that
the Laplace transform of   $\nu^g_t$  is given by
\begin{equation}\label{laplace}
\int_0^{+\infty} 
e^{-sx}\, d\nu_t^g (s)
 =
\int_0^{+\infty} 
e^{-sx} \eta_t^g (s)ds
 =
 e^{-tg(x)}, \;x>0.
\end{equation}
 There is one-to-one correspondence between $g$ and $(\eta^g_t)_{t>0}$, see
\cite{j} p.177, \cite{bf}.

Let recall that all Bernstein functions $g$ can be written by L\'evy-Khintchine formula 
\begin{equation}\label{levy}
g(x)=a+bx+\int_0^{\infty} (1-e^{-\lambda x}) d\nu(\lambda)=a+bx+J(x)
\end{equation}
with $a,b\geq 0$ and $\nu$ a positive measure on $(0,+\infty)$ such that $\int_0^{\infty} \frac{\lambda}{1+\lambda}\, d\nu(\lambda)<\infty$.
The triplet $(a,b,\nu)$ is uniquely determined by $g$, see  Theorem 3.9.4 in  \cite{j} p.174.
We have $a=g(0)=0$ if and only if  $\nu^g_t$ is a probability measure for any $t>0$.
For instance with $a=b=0$, the L\'evy measure associated to $g(x)=x^{\alpha}$ with $ \alpha\in (0,1)$ 
(resp. $g(x)=\ln (1+x)$)
is given by 
$d\nu(\lambda)=\frac{\alpha}{\Gamma(1-\alpha)} \lambda^{-1-\alpha}\,d\lambda$
(resp.  $d\nu(\lambda)=\frac{e^{-\lambda}}{\lambda}\, d\lambda$).

Now, let  $(E_{\lambda})_{\lambda\in  [0,+\infty)}$ be the spectral resolution of the non-negative self-adjoint operator  $A$ and  
$\Psi:[0,+\infty)\longrightarrow \R$   any   measurable function. The operator $ \Psi(A)$ is defined on $L^2(\mu)$ by the 
formula 
  $$
  \Psi(A)f=\int_{0}^{+\infty}   \Psi(\lambda)\, dE(\lambda)f$$
  with domain 
  $$
  {\mathcal D}(   \Psi(A))=\{ f\in L^2: \int_{0}^{+\infty} \vert \Psi(\lambda)\vert^2\, d(E(\lambda)f,f)<\infty\}.
  $$
  We shall set   $( \Psi(A)f,f)=+\infty$ when $f\notin  {\mathcal D}(   \Psi(A))$.
When $\Psi$ has real non-negative values, the operator $ \Psi(A)$ is non-negative and self-adjoint on $L^2(\mu)$.
It defines a symmetric   semigroup of contractions on $L^2$ by the spectral  formula:
  $$
  e^{-t\Psi(A)}f=\int_{0}^{+\infty} e^{-t\Psi(\lambda)}\, dE(\lambda)f,\quad f\in  L^2(\mu).
  $$
 The operator  $e^{-t\Psi(A)}$ will be denoted by   $T^{\Psi}_t$.
 When $\Psi=g$ is a Bernstein function, it can be easily shown that the     semigroup  $(T_t^g)_{t>0}$  satisfies  also the so-called  subordination formula
\begin{equation}\label{subform}
T^{g}_t=\int_0^{+\infty} T_s\, d\nu^g_t(s)=\int_0^{+\infty} T_s \,\eta^g_t(s) ds.
\end{equation}

 Let $t>0$, since $T_t$ is a contraction  on $ L^p$ (resp. positive or  sub-markovian on $L^2$) 
then $T^{g}_t$ is a contraction on $ L^p$ 
(resp. positive or sub-markovian on $L^2$).

Among many examples of Bernstein functions, we are interested at least  in  the following ones:
\begin{enumerate}
\item
The fractional subordinator 
(one-sided $\alpha$-stable process): $g(x)=x^{\alpha}, \; x>0$,
($0<\alpha<1)$. Then
$g(A)=A^{\alpha}$.
\item
The Gamma  subordinator: $g(x)=\log(1+x)$. Then $g(A)=\log (I+A)$ where  $I$ denotes  the identity operator on $L^2(\mu)$.
\item
The generalization of the previous  example:
$g(x)=[\log(1+x^{\alpha})]^{\gamma}$ with $0<\alpha<1$ and $0<\gamma\leq 1$. It gives  $g(A)=[\log (I+A^{\alpha})]^{\gamma}$.
When $\gamma=1$, $g$  is called the geometric $\alpha$-stable subordinator, see \cite{ssv}.
\item
Elementary functions
$g_{\lambda}(x)=1-e^{-\lambda x}$, $\lambda >0$. Then  $g_{\lambda}(A)=
I-T_{\lambda}$.
\end{enumerate}
For a recent study of the    case (ii), see  \cite{ssv}.
 See  also \cite{scsv}  Chap.15 for   a long list of  examples of  Bernstein functions.
 \\

We have  the inclusion ${\mathcal D}(A)\subset {\mathcal D}(g(A))$ and   by formula~\eqref{levy}, for any $f\in{\mathcal D}(A)$,  
\begin{equation}\label{gdea}
g(A)f=af+bAf+\int_0^{\infty} (f-T_{\lambda }f) d\nu(\lambda).
\end{equation}
See \cite{scsv} Example 11.6 (note that  our $A$ is their $-A$).
\\

Throughout all the paper, we shall always assume implicitly that the functions $f$ are in the domain of the operator under consideration. If not, we set $(g(A)f,f)=+\infty$.

We now state the  main results of the paper.


  \begin{ethm} \label{main1}
Let $-A$ be an infinitesimal generator of a
semigroup   as above
and $g$ be a  Bernstein function. If  $-A$ satisfies a  super-Poincar\'e  inequality with rate function $\beta$.  
 
 Then the  infinitesimal generator   $-g(A)$  satisfies  a  $(r_0,r_1)$-super-Poincar\'e  inequality
  with rate function $\beta_{g}(r)=\beta\left( \frac{1}{g^{-1}(1/r)}\right)$, $r\in (r_0,r_1)$ where    $r_0=\frac{1}{g(+\infty)}$ and $ r_1=\frac{1}{g(0^+)}$. 
 \end{ethm}
   
Note that by (\ref{levy}),   $g$ is either strictly increasing    and $g^{-1}$ is defined from $(g(0^+),g(+\infty))$ into
$(0,+\infty)$  or $g$  is constant and $(r_0,r_1)$ is empty. The Bernstein function $g$ is bounded if and only if $b=0$ and $\nu$ is a bounded measure, see \cite{j} p.174. If $g$ is not bounded, 
i.e. $g(+\infty)=+\infty$, then
we have $r_0=0$. If $g(0^+)=a=0$ then $r_1=+\infty$.
But if  $a=g(0^+)>0$ then we have   the  obvious spectral  gap inequality
  $$
  \vert\vert f \vert\vert_2^2\leq \frac{1}{a}(g(A)f,f).
  $$
In other words,
  $$
  \vert\vert f \vert\vert_2^2\leq r(g(A)f,f)+\beta_g(r)
    \vert\vert f \vert\vert_1^2
     \;\; {\rm with}\;\;
     \beta_g(r)= 0,\;  r\in [1/a,+\infty ).
  $$
So, in any case we can   consider that $r_1=+\infty$.
\\

In particular, Theorem \ref{main1}  applies   to several important examples. We provide a short list of  couples $(g(A),\beta_g)$ just below.
\begin{enumerate}
\item
{\em Fractional powers}. If   $g(A)=A^{\alpha}$  with $0<\alpha <1$ then
$
\beta_{\alpha}(r)=\beta( r^{\frac{1}{\alpha}}),\; r>0,
$
(improving  constants given in  \cite{w1}).
\item
{\em  Gamma subordinator}. If $g(A)=\log(I+A)$ then 
$
\beta_{\log}(r)=\beta\left((e^{1/r}-1)^{-1}\right),\; r>0.
$
\item
{\em Generalized geometric stable subordinators}.
If 
$
g(A)=[\log (I+A^{\alpha})]^{\gamma}
$ with  $0<\alpha$, $\gamma\leq 1$,
then 
$$  
\beta_g(r)= \beta\left( \left[ e^{{(\frac{1}{r})}^{\frac{1}{\gamma}}} -1\right]^{\frac{-1}{\;\;\alpha}}\right), \; r>0.
$$
\item
{\em Random walks}.
If
 $g(A)=I-T_{\lambda}=I-e^{-\lambda A}$ with $\lambda>0$ then
$$
\beta^{(\lambda)}
(r)=\beta\left(\frac{\lambda}{\log(1+\frac{1}{r-1})}\right),\; r>1.
$$
The same super-Poincar\'e inequality for the generator $B_{\lambda}=I-T_{\lambda}$ with ${\lambda}>0$ can  also  be  deduced   by a different route,
see (iii) of Proposition \ref{wa} below.
\end{enumerate}

The second main result is  similar   with  Nash-type inequality assumption.

  \begin{ethm} \label{main2}
 Let $-A$ be an infinitesimal generator of a 
semigroup   as above  satisfying  a   Nash-type    inequality with rate function $D$.  Set for $r>0$,  
 $$
 \beta(r)=\sup_{x>0}\left( x-rxD(x)\right),
 $$
and assume that $\beta(r)$ is finite for any $r>0$. Let $g$ be a  Bernstein function. 
\begin{enumerate}
\item
Then
$g(A)$ satisfies a Nash-type inequality of the form
\begin{equation}\label{nashdg}
\vert\vert f\vert \vert_2^2\, D_{g,\beta}(\vert\vert f\vert \vert_2^2)\leq (g(A)f,f), \quad \vert\vert f\vert \vert_1\leq 1,
\end{equation}
where $D_{g,\beta}(x)=\sup_{r\in (r_0,r_1)}\left( \frac{1}{r}-\frac{1}{rx} \beta\left(\frac{1}{g^{-1}(\frac{1}{r})}\right)
\right),\; x>0$  and  $(r_0,r_1)=\left(\frac{1}{g(+\infty)},\frac{1}{g(0^+)}\right)$.
\item
Moreover, assume that $g$ is a bijection from $[0,+\infty)$  to itself  and $\beta:(0,+\infty)\rightarrow (0,+\infty)$ is a decreasing differentiable   bijection.  
Then  we have
\begin{equation}\label{gd}
\vert\vert f\vert \vert_2^2  \,
\sup_{\rho>1} \,( 1-\rho^{-1})  
\; (\,g\circ D\,)( \rho^{-1} \vert\vert f\vert \vert_2^2)
\leq
(g(A)f,f), 
 \quad \vert\vert f\vert \vert_1\leq 1,
\end{equation}
and,  for any $x>0$,
$$
 \sup_{\rho>1} 
( 1-  \rho^{-1}) \,(\,g\circ D\,)(\rho^{-1} x)
\leq D_{g,\beta}(x)
\leq g\circ D(x).
$$
\end{enumerate}
 \end{ethm}

  See Section~\ref{appen2} for the links between $D$ and  $\beta$.
\\

Our approach  simplifies and generalizes  the proofs of  the main results of \cite{bm} and
\cite{w1}. The inequality (\ref{gd}) also clarifies the   constants   obtained in  \cite{bm}  for the fractional powers  $A^{\alpha}$.
 With  the same arguments of proof, we can replace 
$\vert\vert f\vert\vert_1$ in Theorems  \ref{main1} and \ref{main2} by any non-negative functional $\Phi(f)$  satisfying
$\Phi(T_tf)\leq \Phi(f)$, $t>0$.  Our  results can be generalized in the same way in  Hilbert spaces  as in Wang's paper \cite{w1}.  But we shall not  give details.
\\


 \noindent
The paper is organized as follows: 

In Section \ref{intro}, we describe the setting of our study and we state the  main theorems.
In  Section \ref{proo}, we prove the main theorems of Section \ref{intro}. More precisely:
In Section \ref{equivdecay}, we first recall  that the  super-Poincar\'e inequality for $A$ is equivalent to
  the decay   for the corresponding semigroup $(T_t)_{t>0}$.
 Sections \ref{proo1} and \ref{proo2} are devoted respectively to the proof of  Theorem \ref{main1}  
 and  \ref{main2} for $g(A)$ with $g$ a Bernstein function using results of Section \ref{equivdecay}.  

In  Section \ref{sectionultra}, we briefly apply our results to study the eventual ultracontractivity property of subordinated  semigroups.
In Section \ref{exset}, we provide several examples of settings where our results  apply:
1) the  Laplacian on the Euclidean space, 2)
 the Laplace-Beltrami operator on some complete Riemannian manifolds,
3) some hypoelliptic operators on Lie groups
and
4)
 the Ornstein-Uhlenbeck operator on the Euclidean space.

In Section \ref{convex}, we study the same  results
 as in Theorem \ref{main1} and \ref{main2}  concerning  Nash-type and super-Poincar\'e inequalities
 but for   $\Psi(A)$ with $\Psi$  convex.
For these  two   type of inequalities, we use spectral representation of the generator.
  From these results, we deduce      converse implications  of  Theorem \ref{main1} and \ref{main2} by noting that the inverse of a Bernstein function is a convex function. 

In Section \ref{revisit}, we revisit the spectral gap in $L^p(\mu)$ for $g(A)$ using the approach by subordination as in  Theorem \ref{main1}.

In Section \ref{fonclap}, we prove super-Poincar\'e inequalities for $g(\Delta)$
with  $\Delta$ the Laplacian on $\R^n$
 for a larger class of functions than Bernstein functions. 
Our tool is   Fourier analysis as used in the original paper 
by J.Nash \cite{N}. The afferent Nash-type inequalities  can be deduced.
The results are similar to Theorem \ref{main1} and \ref{main2} .

We conclude this paper by   an Appendix  Section \ref{appen2}. 
It concerns the Legendre transform which is underlying  in the equivalence between
Nash-type inequalities and super-Poincar\'e inequalities
through the functions $\beta$ and $D$. 
We weaken the  usual conditions on $\beta$ and $D$ of the  N-functions theory,
see \cite{rr} p.13, more  adapted   to our situation.
We provide examples  used in  Section  \ref{exset}.

\section{Proof of  main Theorems}\label{proo}

\subsection{Super-Poincar\'e  inequality versus  semigroup decay   and Nash-type inequality}\label{equivdecay}
We recall some known result used in the proof of Theorem \ref{main1} and \ref{main2}.
\\

For the proof of Theorem \ref{main1}, we use in a crucial way the following result of F-Y. Wang  
namely the equivalence between Super-Poincar\'e   inequality  for $A$ and  the exponential decay  of the associated  semigroup $(T_t)_{t>0}$,
see \cite{w2} p.230 or \cite{w4} Lemma 3.3.5, see also \cite{w3} p.3 and \cite{w4} p.50 for extended results.
 This is the analogue of the equivalence between the usual exponential decay of a semigroup  and    Poincar\'e inequality.  We recall the proof of this proposition for completeness and, additionally,  we show  (iii) 
  that is the exponential decay of  $T_{t}$ turns out to be the super-Poincar\'e inequality  for the operator
$I-T_{t}$. This operator is  related to  the elementary Bernstein function $g_t(x)=1-e^{-tx}$.

\begin{eprop}\label{wa} 
Let $(T_t)_{t>0}$ be a 
semigroup as in Section \ref{intro} with infinitesimal generator $-A$ and let
$\beta:(0,+\infty)\rightarrow (0,+\infty)$.
  Then the three following inequalities are equivalent:
\begin{enumerate}
\item
For any $f\in {\mathcal D(A)}\cap L^1(\mu)$ and  $r>0$,
 \begin{equation}\label{ttpoin1}
\vert\vert f \vert\vert_2^2 \leq r(Af,f)+\beta(r)\vert\vert f \vert\vert_1^2.
 \end{equation}
\item
For any  $f\in L^2(\mu)\cap L^1(\mu)$, $t>0$ and $r>0$,
   \begin{equation}\label{decaysupoi3}
\vert\vert T_t f\vert\vert_2^2\leq e^{-2t/r}
\vert\vert f \vert\vert_2^2 + (1-e^{-2t/r})\beta(r)\,     \vert\vert f \vert\vert_1^2.
   \end{equation}
\item
For any  $f\in L^2(\mu)\cap L^1(\mu)$, $t>0$ and  $r>1$,
\begin{equation}\label{elementary}
\vert\vert f\vert\vert_2^2 \leq
r((I-T_{t})f,f)+
\beta\left(\frac{t}{ \log (1+\frac{1}{r-1})}\right)
\vert\vert f\vert\vert_1^2.
\end{equation}
\end{enumerate}
\end{eprop}
 
The exponential $ e^{-2t/r}$ in  (\ref{decaysupoi3})  is suitable  to deal with Laplace transforms
(\ref{laplace}) and this  is the key point of our paper.
This allows us  to transfer easily  (\ref{decaysupoi3}) from $A$ to $g(A)$.
\\
\indent
During  the proof, we can notice  that the equivalence between (\ref{ttpoin1})  and
(\ref{decaysupoi3})
holds for any fixed $a=r>0$ and fixed $b=\beta(r)\geq 0$. In particular, if (\ref{ttpoin1})  holds on some interval $(r_0,r_1)$ then  
(\ref{decaysupoi3}) also holds on the same interval $(r_0,r_1)$ and conversely.
\\
\indent
 The inequality (\ref{elementary}) corresponds  exactly  to  a super-Poincar\'e inequality for
 $g_{\lambda}(x)=1-e^{-\lambda x}$ ($t={\lambda}$) for any ${\lambda}>0$.
  The equivalence between (\ref{ttpoin1})  and  (\ref{elementary}) is particularly interesting in terms of relationships between  the super-Poincar\'e profile  for generators $A$ and $B_{\lambda}=I-T_{\lambda}$
  for any fixed ${\lambda}>0$, see  (\ref{profile})  above  for the definition of the profile.
  If $\beta_p(s),s>0$ is the super-Poincar\'e profile of $A$  and $\gamma_p^{(\lambda)}(r), r>1$ is  
  the super-Poincar\'e profile of $B_{\lambda}$  then   they correspond by the formulas
  $$
  \gamma_p^{({\lambda})}(r)=\beta_p\left(\frac{\lambda}{ \log (1+\frac{1}{r-1})}\right), \; r>1,
  $$
  or equivalently
  $$
  \beta_p(s)=\gamma_p^{(\lambda)}\left(1+(e^{\lambda/s}-1)^{-1}\right), \; s>0.
  $$
  For instance in the Euclidean setting, the optimal Nash inequality (\ref{bestnash}) below provides  
  the super-Poincar\'e profile for the Laplacian $\Delta$,
   namely 
  $\beta_p(s)= C_n\,s^{-n/2}$,   for some optimal constant $C_n$. Thus, the super-Poincar\'e profile of $B_{\lambda}$ is explicit and  given by
  $$
  \gamma_p^{(\lambda)}(r)=C_n\,\lambda^{-n/2}
  \left(
  \log \left[1+\frac{1}{r-1}\right]\right)^{n/2}, \; r>1.
  $$
  
  We have the following interpretation in terms of {\em  random walks}.
  For fixed $\lambda>0$, the kernel $h_{\lambda}$ of the operator $T_{\lambda}$ can be seen  as a probability  transition  of a discrete random walk
  $(X_k)_k$  on $\R^n$ given by
  $\dP(X_{k+1}=x,X_k=y)=h_{\lambda}(x-y)=
  \frac{1}{(4\pi {\lambda})^{n/2}} \exp(-\frac{\vert x-y\vert^2}{4{\lambda}})$,
  $
   x,y \in  \R^n
  $
and  $k\in \N\cup\{0\}$.
The  operator  $B_{\lambda}=I-T_{\lambda}$    is the generator of the continuous-time Markov semigroup
  $Q^{(\lambda)}_t=e^{-tB_{\lambda}}=e^{-t}\sum_{k\geq 0} \frac{t^k}{k!} T_{k{\lambda}}
  $
 obtained by convolution  with  the following probability  transition
   $$ q^{(\lambda)}_t 
   =
  e^{-t}\sum_{k\geq 0} \frac{t^k}{k!}  
  h_{k\lambda} 
  $$
  where $h_0=\delta_0$ ($\delta_a$ is the Dirac mass at $a\in [0,+\infty$)).
This semigroup $Q^{(\lambda)}_t$ is subordinated to the heat semigroup  $e^{-t\Delta}$ 
by the  Poisson semigroup  with jumps of size  ${\lambda}$ defined on $[0,+\infty)$ by
$\nu_t =
\sum_{k\geq 0} \frac{t^k}{k!} e^{-t}\, \delta_{k{\lambda}}
$
in  (\ref{subform}), see  \cite{j} p.180.
\\
   
 
{\bf Proof:}    
{\em Equivalence between (i) and (ii)}. 
 Let      $H(t)=e^{2t/r}\vert\vert T_t f\vert\vert_2^2$ for $t >0$, fixed $r>0$
and
   $f\in {\mathcal D(A)}\cap L^1(\mu)$.     
We have
$$
H(t)-H(0)=\int_0^tH'(u)du=
\int_0^t
2\, e^{2u/r}\left( \frac{1}{r} \vert\vert T_u f\vert\vert_2^2 -(AT_uf,T_uf)\right)du.
 $$
By applying (\ref{ttpoin1})  to $T_uf$ and   since  $T_u$ is a contraction on $L^1(\mu)$, we 
deduce
$$
H(t)-H(0)
\leq 
 \frac{2}{r} \beta(r) \vert\vert f\vert\vert_1^2
\left( 
\int_0^t
e^{2u/r}
 \,du
 \right).
 $$
 This proves   (\ref{decaysupoi3}) for $f\in {\mathcal D(A)}\cap L^1(\mu)$. For the general case,
 let $f\in L^2(\mu)\cap L^1(\mu)$ then  there exists $f_u$ $(u>0)$  such that $f_u\in {\mathcal D(A)}$  and $f_u$ converges to $f$ in $ L^1(\mu)$ and $ L^2(\mu)$
 as $u\rightarrow 0^+$
 (e.g. $f_u=\frac{1}{u}\int_0^u T_sf\,ds$).

Conversely, let $r>0$ be fixed and $f\in {\mathcal D(A)}\cap L^1(\mu)$. The inequality  (\ref{decaysupoi3}) can be rewritten as
 $$
\frac{ \vert\vert T_t f\vert\vert_2^2 -  \vert\vert f\vert\vert_2^2}{2t}
\leq
\left(  \frac{e^{-2t/r} -1}{2t}\right) 
\vert\vert f \vert\vert_2^2 +
\beta(r)\, \vert\vert f \vert\vert_1^2
    \left(\frac{1-e^{-2t/r}}{2t}\right).
 $$
 We conclude   (\ref{ttpoin1}) by taking the limit as $t$ goes to $0$.

{\em Equivalence between (ii) and (iii)}.
Assume that  (ii) holds, i.e.
$$
(T_tf,f)= \vert\vert T_{t/2} f\vert\vert_2^2\leq e^{-t/r}
\vert\vert f \vert\vert_2^2 + ({1-e^{-t/r}} )  \beta(r)\,
 \vert\vert f \vert\vert_1^2
 $$
  for any $t$, $r>0$.
It is equivalent to   
    $$
   ({1-e^{-t/r}} ) \vert\vert f \vert\vert_2^2 
    \leq
   (f-T_{t} f,f)
    +
     ({1-e^{-t/r}} )  \beta(r)\,
   \vert\vert f \vert\vert_1^2.
 $$
 Let $g(x)=g_t(x)=1-e^{-tx}$. The last inequality reads as
 $$
\vert\vert f \vert\vert_2^2 
    \leq
    \frac{1}{ g_t(1/r)}
   (g_t(A) f,f)
    + \beta(r)\,
   \vert\vert f \vert\vert_1^2.
 $$
Fix $t>0$.
 Let $\rho>1$ and choose $r>0$ such that $\rho=  \frac{1}{ g_t(1/r)}$, i.e.
 $r=\frac{1}{ g^{-1}_t(1/\rho)}=\frac{t}{\log (1+\frac{1}{\rho-1})}$.
This  yields  
 $(1,\infty)$-super-Poincar\'e  (\ref{elementary}) for the operator $I-T_{t}$
as expected.
The converse is clear.
   \rule{1.3mm}{2mm}
 \\


 Now we recall that super-Poincar\'e and Nash-type inequalities  are essentially equivalent  under natural conditions on $\beta$ in (\ref{poini}) and on $D$ in (\ref{nashi}). 
 This result is more or less well known   but 
 we formulate  the relations between  $\beta$  and $D$ implicitly in terms of Legendre transforms,  
see    Appendix  Section \ref{appen2} for a  detailed discussion.  
\begin{eprop}\label{nashpoin}
Let $A$ be a non-negative  symmetric  operator on $L^2(\mu)$.
\begin{enumerate}
\item
 Assume $A$ satisfies  a super-Poincar\'e inequality with rate function $\beta$
then it satisfies a Nash-type inequality with rate function
$$
D(x)=\sup_{t>0}\left(t-\frac{t \beta(1/t) }{x}\right)\in (-\infty,+\infty],\; x>0.
$$
 The function $D$
is non-decreasing, finite on the set $(0,\sup{\mathcal G})$ where 
${\mathcal G}=\{ \vert\vert f \vert\vert_2^2, f\in {\mathcal D}(A)\cap L^1(\mu), \vert\vert f \vert\vert_1\leq 1\}$
and 
$D(+\infty)=+\infty$.
\item
Conversely, suppose  Nash-type inequality holds true for  $A$ and a rate function   $D: [0,+\infty)\rightarrow \R$. 
Set 
\begin{equation}\label{betaf}
\beta(r)=\sup_{x>0} \left(x-rxD(x)\right)\in (-\infty,+\infty].
\end{equation}
 Then $A$ satisfies a super-Poincar\'e  inequality with rate function $\beta$.
 \end{enumerate}
\end{eprop}
This proposition is   in the spirit of Theorem 3.1 and Section 5 of \cite{w2}, see also Proposition 3.3.16 of \cite{w4}.
A  much closer formulation  can be found in  \cite{bim}, Proposition 2.1.
\\

Note that if ${\mathcal G}$ of   Proposition \ref{nashpoin} is unbounded above as a subset of $\R$ then   $D$ is finite in $(0,\infty)$. 
For many  examples, $\beta$ is a non-negative  non-increasing function and satisfies $\beta(0^+)=+\infty$, $\lim_{t\rightarrow 0^+} t\beta(1/t)=0$ which implies that $D(x)$ is finite for any $x>0$, non-negative, continuous and non-decreasing, see Appendix Section \ref{appen2} for details.
Note also that we can always consider in (\ref{nashi}) that $D$ is non-negative by  replacing  $D$ by $D_+=\sup(0,D)$ 
since  $A $ is a non-negative operator.

\begin{erem}
\begin{enumerate}
\item Usually Nash-type inequality is written in the following form
$$
\Theta(\vert\vert f \vert\vert_2^2) \leq (Af,f), \quad \vert\vert f \vert\vert_1\leq 1.
$$
 But the equivalent expression \eqref{nashi} is more appropriate to deal with the Bernstein functions $g(x)=x^{\alpha}$, $\alpha \in (0,1)$ as shown in  \cite{bm} 
and more generally for any  Bernstein functions by Theorem \ref{main2} above.
\item
In the second statement, the assumption on the functional $(Af,f)$ can certainly be relaxed  because  only the existence of the function $\beta$ (which depends on  this functional) is   crucial for  the proof.
\end{enumerate}
\end{erem}
{\bf Proof:}
  (i) 
  Assume that 
a super-Poincar\'e inequality   holds true.
For any  $\vert\vert f \vert\vert_1\leq 1$, $f\neq 0$, $ f\in {\mathcal D}(A)$  and  any $r>0$, we easily deduce
$$
  \vert\vert f \vert\vert_2^2\left( \frac{1}{r} -\frac{\beta(r)}{r \vert\vert f \vert\vert_2^2}\right)\leq (Af,f).
$$
Taking the supremum over $r>0$,  we get a Nash-type  inequality with rate function $D$.    
  Note that 
$D$ is automatically finite on the subset  ${\mathcal G}\setminus \{0\}$ of $\R$ since in that case $(Af,f)<+\infty$
for  any $ f\in {\mathcal D}(A)$. 
On one hand, the set ${\mathcal G}$ is not empty since it contains $0$. On the other hand, if $f\in {\mathcal G}$ and 
$\lambda\in (0,1)$ then $\lambda f\in {\mathcal G}$. Hence $(0, \sup {\mathcal G})\subset  {\mathcal G}$
and
 $D$ is finite on $(0,\sup{\mathcal G})$.
 It is easily proved that $D$  is non-decreasing using   the fact that ${\beta}>0$.
  Moreover, for   any $x$, $t>0$, we have
$$
D(x)\geq t-\frac{t\beta(1/t)}{x}.
$$
Therefore, 
$\displaystyle\liminf_{x\rightarrow + \infty}D(x)\geq t$ for any  $t>0$.
It  implies
$\displaystyle\lim_{x\rightarrow +\infty}D(x)=+\infty$.
\\
 (ii)  
  By definition of $\beta$, one has for any $x$, $r>0$,
$$
\frac{x}{r}-\frac{{\beta}(r)}{r}\leq  xD(x).
$$
Let $x=\vert\vert f  \vert\vert_2^2$ (with $f\neq 0$).
So, for fixed $r>0$,
$$
\frac{\vert\vert f  \vert\vert_2^2}{r}- \frac{1}{r}{\beta}(r)\leq \vert\vert f  \vert\vert_2^2\,D(\vert\vert f  \vert\vert_2^2).
$$
By the Nash-type inequality, the last term is bounded by  $(Af,f)$ when $\vert\vert f  \vert\vert_1\leq 1$. Hence,
$$
\frac{\vert\vert f  \vert\vert_2^2}{r}- \frac{1}{r}{\beta}(r)\leq  (Af,f)
$$
which is  super-Poincar\'e inequality. 
This concludes the proof of the proposition.
 \rule{1.3mm}{2mm}
\\

Note that   when $\beta(r)=+\infty $  for some $r>0$ then  super-Poincar\'e inequality  is satisfied and the proof is also valid. 
For properties of $\beta$ defined by (\ref{betaf}), see Appendix Section \ref{appen2}.

\subsection{Proof of  Theorem \ref{main1}}\label{proo1}
We suppose that  $g$ is non-constant (if not  there is nothing to prove since $(r_0,r_1)=\emptyset$).
Assume that a super-Poincar\'e inequality holds true with rate function $\beta$.  
By Proposition~\ref{wa}, the inequality (\ref{decaysupoi3}) is satisfied.
Now by   symmetry and   semigroup property of $(T_t)_{t>0}$, this inequality  (\ref{decaysupoi3}) can be written as
$$
\vert\vert T_tf \vert\vert_2^2=( T_{2t} f,f) \leq e^{-2t/r}
\vert\vert f \vert\vert_2^2 +\beta(r)\, \vert\vert f \vert\vert_1^2\,
    (1-e^{-2t/r}), \; t>0.
$$
Let  $s>0$ and set $t=s/2$. We deduce for any $r>0$ and 
$s>0$,
$$
( T_{s} f,f)\leq e^{-s/r}
\vert\vert f \vert\vert_2^2 +\beta(r)\, \vert\vert f \vert\vert_1^2\,
    (1-e^{-s/r}).
$$
By   the subordination formula (\ref{subform}) and Fubini, we get for any $t$, $r>0$, $f\in L^1(\mu)\cap L^2(\mu)$ and  any Bernstein function $g$,
$$
(T_t^g f,f)
=
  \int_0^{+\infty} 
 (T_sf,f) \, d\nu^g_t(s)
\leq
\left( \int_0^{+\infty}  e^{-s/r}  \, d\nu^g_t(s)\right) 
 \vert\vert  f\vert\vert_2^2 
 $$
 $$
 \;+\; \beta(r)\,  \vert\vert  f\vert\vert_1^2
  \left(   \int_0^{+\infty} (1-e^{-s/r}) \, d\nu^g_t(s) \right).
$$
By the Laplace transform of the sub-probability $\nu^g_t$, we get
$$
(T_t^g f,f)
  \leq e^{-tg\left(1/r\right) }
\vert\vert   f\vert\vert_2^2 +\beta(r)\,\vert\vert   f\vert\vert_1^2
    (1-e^{-tg\left(1/r\right)}), \quad  t, \,r>0. 
$$
Changing $t$ by $2t$ and using  symmetry and semigroup properties of $(T^g_t)$, we obtain 
$$
 \vert\vert T_t^g f\vert\vert_2^2  
    \leq e^{-2tg\left(1/r\right) }
\vert\vert   f\vert\vert_2^2 +\beta(r)\,\vert\vert   f\vert\vert_1^2
    (1-e^{-2tg\left(1/r\right)}), \quad  t, \,r>0. 
$$
Now let $\rho\in (r_0,r_1) :=(\frac{1}{g(+\infty)},\frac{1}{g(0^+)})$. 
Since  $r\longrightarrow \frac{1}{g\left(1/r\right)}$ is a bijection from $(0,+\infty)$ onto $(r_0,r_1)$, there exists a (unique)  $r>0$ such that ${\rho} =\frac{1}{g\left(1/r\right)}$, i.e.
$r=\frac{1}{g^{-1}\left(1/\rho\right)}$, 
and
$$
\vert\vert T_t^g f\vert\vert_2^2
  \leq e^{-2t/\rho}
\vert\vert   f\vert\vert_2^2 +\beta \left(\frac{1}{g^{-1}\left(1/\rho\right)}\right)\, \vert\vert   f\vert\vert_1^2
    (1-e^{-2t/\rho}). 
   $$
   We conclude  by 
applying   (ii) $\Rightarrow$ (i)  of Proposition   \ref{wa} with $g(A)$.
Theorem~\ref{main1} is proved.
\rule{1.3mm}{2mm}

Note that we do not need the existence of the density   of the measures $\nu^g_t$  
 nor additional properties of the function $\beta$.



\subsection{Proof of  Theorem \ref{main2}}\label{proo2}
{\em  Proof of (i)}.
We assume that $A$ satisfies Nash-type inequality.
By      (ii) of Proposition \ref{nashpoin}  and   the  definition  of $\beta$, we get super-Poincar\'e
inequality: for any $f\in {\mathcal D}(A)\cap L^1(\mu)$,
$$
 \vert\vert f \vert\vert_2^2   
 \leq 
 r(Af,f)+ \beta(r)\vert\vert f \vert\vert_1^2,\, r>0,
 $$
We now apply   Theorem \ref{main1} and deduce super-Poincar\'e  inequality for $g(A)$, i.e. 
$$
\vert\vert f \vert\vert_2^2 \leq r(g(A)f,f)+\beta_g(r)\vert\vert f \vert\vert_1^2, 
$$
   with $\beta_{g}(r)=\beta\left( \frac{1}{ g^{-1}(1/r)}\right)$, 
   $r\in (r_0,r_1)=(\frac{1}{g(+\infty)},\frac{1}{g(0^+)})$. 
Now, we conclude by applying    (i) of  Proposition \ref{nashpoin}.
\\

\noindent
{\em  Proof of  (ii)}. 
From the next lemma which compares $D_{g,\beta}$ and $g\circ  D$ and (i), we immediately deduce the inequality (\ref{gd}).

\begin{flem} 
Let $g:(0,+\infty)\rightarrow (0,+\infty)$ be a bijective  continuous increasing concave function
(e.g. bijective Bernstein function), $D$ and  $D_{g,\beta}$  defined as in Theorem~\ref{main2} with $\beta$  a decreasing  differentiable bijection  from  $(0,+\infty)$ to itself.
 Then for any $ x>0$ and $\rho>1$,
$$
 (1- {\rho}^{-1})
 (g\circ D)({\rho}^{-1} x )
\leq D_{g,\beta}(x)\leq g\circ D(x).
$$
\end{flem}

{\bf Proof:}
To simplify our discussion,  we set $V(t):=\beta(1/t) $, $V_g(t):=\beta_g(1/t)=\beta(\frac{1}{g^{-1}(t)})$ for $t>0$.
So,
$
D(x)=\sup_{t>0}\left( t-\frac{t}{x}V(t)\right)
$
and 
$
D_{g,\beta}(x)=\sup_{t>0}\left( t-\frac{t}{x}V_g(t)\right)$ for $x>0$. 
As a consequence of the assumptions $\beta(0^+)=+\infty$, we have that $D$ and $D_{g,\beta}$ are well defined and finite  on $(0,+\infty)$.

 Let $u>0$. Since $g$ is a bijection from $(0,+\infty)$ to itself, there exists a unique $t>0$ such that 
 $\frac{1}{u}=\frac{1}{g^{-1}(t)}$,
 i.e 
 $t= g({u})$. Thus $D_{g,\beta}$ can be written as 
 $$
 D_{g,\beta}(x)=\sup_{u>0} \; g(u)\left(1-\frac{V(u)}{x}\right).
 $$
Since $D=D_{id}$ and by continuity of $g$, we get
 $$
 g\circ D(x)=  
 \sup_{\{u>0: V(u)\leq x\}}\,  g\left( u\left[1- \frac{V(u)}{x}\right]\right) .
 $$
Let $a=1- V(u)/x$. 
Since $V\geq 0$, it is sufficient  to consider  $a\in (0,1)$. 
By concavity of $g$ and $g(0)=0$, we have $ag(u)=ag(u)+(1-a)g(0)
\leq g(au)$.
 Therefore, we conclude    $D_{g,\beta}(x)\leq g\circ D(x)$ for any $x>0$.

Now, we prove the lower bound on $D_{g,\beta}$.
From the  definition, we have   for any $x$, $u>0$,
 $$
 D_{g,\beta}(x)\geq  g(u)\left(1-\frac{V(u)}{x}\right).
 $$
By the assumptions on  $\beta$, the function $V$ is
a differentiable   increasing  bijection from $(0,+\infty)$ to itself. 
Fix $x>0$. For $\rho>1$, we set  $u=V^{-1}( {\rho}^{-1} x )$.
It yields
 $$
  D_{g,\beta}(x)\geq (1- {\rho}^{-1})
  \,g(V^{-1}( {\rho}^{-1} x) ). 
  $$

Fix $y>0$. The supremum defining
$D(y)=t_0- t_0 y^{-1}V(t_0)$ exists  and it is attained at some  point $t_0>0$  which  is characterized by $1-\frac{1}{y} V(t_0)-\frac{t_0}{y}V'(t_0)=0$.
It implies that
$y= V(t_0)+t_0V'(t_0)\geq V(t_0)$ because $V'\geq 0$. Finally, we get 
$t_0\leq V^{-1}(y)$.
Since $V\geq 0$, we deduce 
$D(y)\leq t_0\leq V^{-1}(y)$. Thus $g(D(y))\leq g(V^{-1}(y))$ for any $y>0$. Now set $y={\rho}^{-1} x $ and obtain the expected lower bound
$$
 D_{g,\beta}(x)
  \geq 
 (1-{\rho}^{-1})
 (g\circ D)({\rho}^{-1} x).
$$
The proof is complete. 
\rule{1.3mm}{2mm}

\section{Application to ultracontractivity of   subordinated semigroups}\label{sectionultra}
Recall  that  a  symmetric semigroup $ (T_t)_{t>0}$ of contraction on $L^2(\mu)$ and $L^1(\mu)$  is    ultracontractive if for any $t>0$,
\begin{equation}\label{ult}
\vert\vert T_tf\vert\vert_2\leq b(t) \vert\vert f\vert\vert_1
\end{equation}
 for some  non-increasing  function $b:(0,+\infty)\rightarrow (0,+\infty)$ with $b(0^+)=+\infty$,
 see \cite{d}.
 Ultracontractivity  implies  super-Poincar\'e~\eqref{poini} with 
$\beta(r)=b^2(r/2)$.  
Indeed,  since  $s\rightarrow (AT_sf,T_sf)$ is non-increasing, we get for any   $r>0$,
$$
\vert\vert f\vert\vert_2^2
-b^2(r/2)\vert\vert f\vert\vert _1^2
\leq
\vert\vert f\vert\vert_2^2-\vert\vert T_{r/2} f\vert\vert_2^2=
(f-T_rf,f)=\int_0^r (AT_sf,T_sf)\,ds \leq r (A f,f).
$$
Which is the desired  inequality.

By interpolation and duality, the property of ultracontractivity is   equivalent 
 to
 \begin{equation}\label{ultinf}
\vert\vert T_tf\vert\vert_{\infty}\leq a(t) \vert\vert f\vert\vert_1,\quad t>0,
\end{equation}
 for some  non-increasing  function $a:(0,+\infty)\rightarrow (0,+\infty)$ with $a(0^+)=+\infty$.  More precisely, from (\ref{ult})
 we get
 $a(t)\leq b^2(t/2)$ and from (\ref{ultinf}) we obtain $b(t)\leq \sqrt{a(t)}$.

If $b_g(t):= \int_0^{+\infty} 
b(s)  \eta_t^{g}(s)   \,ds<+\infty$ for any $t>0$ then the semigroup $(T_t^g)$ is ultracontractive  since 
$$ 
  \vert\vert  e^{-tg(A)} f \vert\vert_2
  \leq
\int_0^{+\infty} 
 \eta_t^{g}(s) \vert \vert T_sf 
\vert \vert _2 \,ds 
\leq
\left( 
 \int_0^{+\infty} 
 b(s)  \eta_t^{g}(s)
\,ds \right) \vert\vert f\vert\vert_1.
 $$
 But unfortunately, to check this condition is rather hard because  the densities  $ \eta_t^{g}$ are not   well known apart from the case of $g(x)=\sqrt{x}$, see \cite{j} p.181. 
 A way to overcome this difficulty is 
  by considering  Nash-type inequalities.
For that purpose, we recall a result due to T.Coulhon. The author    deduces  ultracontractivity bounds  from Nash-type inequality under some integrability condition, see \cite{c} and  also \cite{m1}. For   applications we have in mind, 
 we restrict his result to our setting. 
 
 \begin{ethm}\label{coulhon}
Let $(T_t)_{t>0}$ be  a  semigroup  as in Section \ref{intro} with 
 infinitesimal  generator $-A$.
 Assume that  there exists a non-decreasing function $\Theta:(0,+\infty)\rightarrow (0,+\infty)$    satisfying the following  Nash-type inequality
 $$
 \Theta(\vert\vert f\vert\vert_2^2)\leq (Af,f),\quad \vert\vert f\vert\vert_1\leq 1.
 $$
If
 $\int^{\infty} \frac{dx}{ \Theta(x)}<+\infty$
 then $(T_t)_{t>0}$ is ultracontractive and for any $t>0$,
 $$
 \vert\vert T_tf\vert\vert_{\infty}\leq a(t) \vert\vert f\vert\vert_1,\quad t>0,
 $$
where  $a(t)$ is  the inverse  of the function  $s\rightarrow \int_s^{\infty} \frac{dx}{ \Theta(x)}$.
 \end{ethm}
 We apply this  result to the eventual  ultracontractivity  of the subordinated semigroup $(T^g_t)$
 and  give  a sufficient condition on $D_{g,\beta}$  of (\ref{nashdg}) to get ultracontractivity from a Nash-type inequality satisfied by $A$.

\begin{ecor}\label{corult}
Under the assumptions and notations of Theorem \ref{main2},
let's assume that $\int^{\infty} \frac{dx}{ xD_{g,\beta}(x)}<+\infty$.
Then $(T_t^g)$ is ultracontractive and for any $t>0$,
$$
 \vert\vert T_t^gf\vert\vert_{\infty}\leq a_g(t) \vert\vert f\vert\vert_1
 $$
 where  $a_g$  is  the inverse function of $s\rightarrow \int_s^{\infty}  \frac{dx}{ xD_{g,\beta}(x)}$.
\end{ecor}

{\bf Proof}: Apply Theorems \ref{main2} and \ref{coulhon}.
\\

It is clear that $(T_t)_{t>0}$ can be ultracontractive but not  $(T^g_t)$ for some Bernstein functions $g$. For instance, let 
 $A=\Delta$ be  the usual Laplacian on   $\R^n$.   Then $A$ satisfies a  Nash-type  inequality with rate function $D(x)=cx^{2/n}$.
Let  $g(r)=\log(1+r)$. Then $\int^{\infty} \frac{dx}{ xD_{g,\beta}(x)}\geq 
\int^{\infty} \frac{dx}{ x\ln (1+cx^{2/n})}
=+\infty$.
This is  obtained  from   the inequality $ g\circ D(x)\geq D_{g,\beta}(x)$ of Theorem~\ref{main2} (ii). 
Now by a direct computation, we can   show that $(T^g_t)$ is not ultracontractive for small $t>0$,  
see (\ref{ultboundgeom}) below for details.
\\

{\em Applications to heat kernel bounds:}
Ultracontractivity  insures    the existence and uniform bounds of the heat kernel under some assumptions  on $X$.  
 For instance,  if 
 $X$ is  a locally compact separable metric space  with a Radon measure $\mu$  with full support
then  ultracontractivity  of $(T^g_t)$ implies  existence of the heat kernel $k^g_t$ with respect to the measure $\mu$, 
i.e.
$$
  T^g_t f(x)=\int_X  f(y)k^g_t(x,y)\,d\mu(y).
  $$
Furthermore, the kernel  satisfies  the uniform bound, 
  $$
  {\rm essup}_{x,y\in X} \, k^g_t(x,y)\leq a_g(t), \; t>0.
  $$
See the   recent paper \cite{gh} (Lemma 3.7) for a detailed exposition
on the  existence of the heat kernel. See also \cite{d} Chap.2.

\section{Examples of Settings}\label{exset}
 
 Here, we give  some examples where our results can be applied.
 
 \subsection{The Euclidean space}\label{eucli}
 
 Let  $\Delta=-\sum_{i=1}^n\frac{\partial^2}{\partial x_i^2}$ be the  usual   Laplacian  on $\R^n$. 
 The profile of the super-Poincar\'e  inequality can be deduced from 
 the optimal  Nash inequality obtained in \cite{cl}.
 Let $N_n$ be the best constant in Nash inequality,
 $$
\frac{1}{N_n}  \vert\vert f\vert\vert_2^{2+4/n}\leq (\Delta f,f),\quad \vert\vert f\vert\vert_1\leq  1.
 $$
 By Proposition \ref{nashpoin}, this is equivalent to  the following super-Poincar\'e inequality,
\begin{equation}\label{bestnash}
 \vert\vert f\vert\vert_2^2\leq r(\Delta f,f)+ C_n\,r^{-n/2} \vert\vert f\vert\vert_1^2,\; r>0,
\end{equation}
with 
\begin{equation}\label{bestct}
C_n=\frac{2(nN_n)^{n/2}}{ (n+2)^{1+n/2}}.
\end{equation}
Thus  (\ref{poini}) is satisfied with the super-Poincar\'e profile  $\beta_p(r)= C_n\,r^{-n/2}$ and
(\ref{nashi}) with $D(x)=\frac{1}{N_n} x^{2/n}$.

To simplify the presentation of our results, we shall assume that the Bernstein function $g$   is a bijection from $(0,+\infty)$ to itself.
By applying (ii) of Theorem  \ref{main2}, we get for any $\rho>1$,
 $$
  \frac{1}{2}( 1- {\rho}^{-1})  \vert\vert f\vert\vert_2^2\,g\left( 2{N_n}^{-1}\rho^{-2/n} \vert\vert f\vert\vert_2^{4/n}\right)
   \leq (g(\Delta) f,f),\quad \vert\vert f\vert\vert_1\leq  1.
 $$  

 {\em Examples of Bernstein functions}.
 \begin{enumerate}
 \item
Let  $g(x)=x^{\alpha}$. We obtain for the fractional power of the Laplacian $\Delta^{\alpha},\,0<\alpha<1$,
 $$
  \frac{1}{2 ({N_n} \rho )^{\alpha} }\,( 1-{\rho}^{-1})\,  \vert\vert f\vert\vert_2^{2+4{\alpha}/n}
   \leq (\Delta^{\alpha} f,f),\quad \vert\vert f\vert\vert_1\leq  1.
 $$
 By optimizing over $\rho >1$,
 we deduce 
$$
L_{n,\alpha}\,
\vert\vert f\vert\vert_2^{2+4{\alpha}/n}
\leq (\Delta^{\alpha} f,f),\quad \vert\vert f\vert\vert_1\leq  1,
 $$
 with $L_{n,\alpha}= 2^{\alpha-1} {N_n}^{-\alpha}  
 \frac{
 n({ 2\alpha})^{2\alpha/n}
 }
 { 
 (2\alpha+n)^{1+2{\alpha}/n}
 }
 $.

See  \cite{vsc}
for such a result   in the   setting
of sub-markovian symmetric semigroups  but with  a different approach. Note that the  constant
$
L_{n,\alpha}
$
is explicit but  probably not optimal.
Indeed, we get a better constant if we apply Theorem \ref{main1} with $g(x)=x^{\alpha}$:
 $$
 \vert\vert f\vert\vert_2^2\leq r(\Delta^{\alpha} f,f)+ C_n\,
 2^{\frac{n}{2}(\frac{1}{\alpha} -1)}
 r^{-\frac{n}{2\alpha}} \vert\vert f\vert\vert_1^2,\; r>0,
$$
with $C_n$ as above.
By applying (i) of Proposition \ref{nashpoin}, we get 
$$
    K_{n,\alpha}\,
     \vert\vert f\vert\vert_2^{2+4{\alpha}/n}
   \leq (\Delta^{\alpha} f,f),\quad \vert\vert f\vert\vert_1\leq  1.
 $$
 with 
 $$
 K_{n,\alpha}=
 \left(\frac{n}{n+2\alpha}\right)2^{\alpha -1}\left(\left(\frac{n}{2\alpha}+1\right) C_n\right)^{\frac{-2\alpha}{\; n}}.
$$
By the relationships connecting
  $N_n$ and $C_n$, we have that $L_{n,\alpha}< K_{n,\alpha}$  (equivalent   to the trivial inequality
$n\, 2^{\frac{2}{n}}< (n+2)^{\frac{2}{n}+1}$).
 We postpone to Section \ref{fonclap}  the study of  super-Poincar\'e inequalities for 
 a larger class of functions of the Laplacian    using Fourier analysis tools.
 But with this approach, the best constants are lost.
 
 \item
Let  $g(x)=\log(1+x)$.
 The geometrically stable operator $\log(I+\Delta)$ satisfies 
 $$
 \frac{1}{2}( 1-{\rho}^{-1})  \vert\vert f\vert\vert_2^2\,\log\left( 1+{N_n}^{-1}\rho^{-1} \vert\vert f\vert\vert_2^{4/n}\right)
   \leq (\log(I+\Delta) f,f),\quad \vert\vert f\vert\vert_1\leq  1.
 $$
To estimate
 $D_{g,\beta}$ with $g(x)=\log (1+x)$ is not a pleasant task. So,  we prefer to state this explicit inequality for each  parameter $\rho>1$.
\end{enumerate}
Note that, in general, the eventual ultracontractivity can be proved for $e^{-tg(\Delta)}$ directly by the formula,
$$
\vert\vert e^{-tg(\Delta)}\vert\vert_{1\rightarrow 2}^2=\frac{1}{(2\pi )^{n}}
\int_{\R^n} e^{-2tg(\vert y\vert^2)}\, dy.
$$
Applied to $g(x)=\log(1+x)$,   this  leads us to
\begin{equation}\label{ultboundgeom}
\vert\vert e^{-t\log(I+\Delta)}\vert\vert_{1\rightarrow 2}^2
=
\frac{\vert S_{n-1}\vert}{(2\pi )^{n}}
\int_0^{\infty} (1+r^2)^{-2t}\,r^{n-1}\,dr <+\infty \quad {\rm iff} \quad t>\frac{n}{4}.
\end{equation}
Thus this semigroup is not ultracontractive for $0<t\leq n/4$.
But  note that it satisfies  super-Poincar\'e and Nash-type inequalities.
 
 \subsection{The Riemannian setting}\label{riem}
  
 The following example is taken from \cite{w2} Cor.2.5.
 Let $M$ be a  connected  complete Riemannian manifold with Ricci curvature
bounded from below. Assume that the boundary  $\partial M$ is convex or empty.
For $V\in C^1(M)$, we assume $Z=\int_M e^{-V(x)}\,dx$ is finite and define the probability  measure $\mu$ by
$d\mu(x)=Z^{-1}e^{-V(x)}\,dx$ where  $dx$ is the Riemannian volume measure.
Let  $A=\Delta+\nabla V$, $A$ is (essentially) self-adjoint on $L^2(\mu)$
(with Neumann boundary condition whenever $\partial M$ is nonempty).
 Set $\rho(x)=\rho(x,o)$ the Riemannian distance function to a fixed point $o\in M$.
Consider  $  V=-\alpha\rho^{\delta}$, $\alpha>0$ and $\delta >1$ then 
  super-Poincar\'e  (\ref{poini}) holds true with 
\begin{equation}\label{exporiem}
  \beta(r)=\exp\left[ c\,(1+r^{-\lambda})\right]
\end{equation}
 with $\lambda=\delta/[2(\delta-1)]$ and  some constant $c >0$.  
  Moreover, a super-Poincar\'e holds if $V=-\exp[\alpha \rho],\alpha>0$ with $\lambda=1/2$ and the rate function $\beta$ given in~\eqref{exporiem}.  Theorem \ref{main1} implies that $g(A)$ satisfies a  super-Poincar\'e  with rate function  $\beta_{g}(r)$ given in  Theorem~\ref{main1}.
 Nash-type inequalities  can be deduced  for $g(A)$ from super-Poincar\'e by   Theorem \ref{main2}. 
Our  result generalize the particular case of the fractional powers $g(x)=x^{\alpha}$, $0<\alpha<1$ treated in~\cite{w1}.
 
 \subsection{The hypoelliptic setting}\label{hypo}
  
  Here, we consider sub-laplacians on Lie groups of polynomial growth. Let $G$ be a connected Lie group of polynomial  growth of index $D$ and 
  $(X_1,X_2,...,X_m)$ be a system of left-invariant vector fields satisfying H\"ormander's condition with local dimension $d$. We assume $d\leq D$.
The   sub-laplacian
 $L=-\sum_{i=1}^m X_i^2$  generates a semigroup $e^{-tL}$ with density kernel $p_t$ satisfying for all $n$ satsifying  $d\leq n\leq D$,
 $$
 \sup_{x,y\in G} p_t(x,y)=\vert\vert e^{-tL}\vert\vert_{1\rightarrow \infty}\leq
 \frac{c_1}{t^{n/2}}.
 $$
 
    Hence,  a super-Poincar\'e inequality holds true
 $$
 \vert\vert f\vert\vert _2^2\leq t(Lf,f)+ \frac{c_0}{t^{n/2}} \vert\vert f\vert\vert _1^2, \quad t>0.
 $$
and our results  applies  to $g(L)$ for any Bernstein function $g$  with $\beta(t)= c_0 t^{-n/2}$, see \cite{vsc}.
We now discuss in this context the four examples of Bernstein function $g$  introduced below  Theorem \ref{main1}.
We provide  asymptotic behaviors of $\beta_g(r)$ when $r$ tends to $0$ and $r$  tends  to $+\infty$.
\begin{enumerate}
\item
If $g(x)=x^{\alpha}, 0<\alpha \leq 1$ then $\beta_g(r)=\frac{c_0}{r^{n/2\alpha}}, \;r>0 $.
\item
If  $g(x)=\ln (1+x) $ then $\beta_g(r)= c_0 \,(e^{1/r}-1)^{n/2},\; r>0$.
$$
 \beta_g(r) \sim 
 \begin{cases}
 c_0 \;e^{n/2r} \;\;{\rm as}\; \;\, r\rightarrow 0^+ ,\\
c_0 \; \frac{1}{r^{n/2}}   \;\;{\rm as}\;\; r\rightarrow +\infty.
\end{cases}
$$
\item
If  $g(x)=\left[\ln (1+x^{\alpha})\right]^{\gamma}, \; 0<\alpha$, $\gamma\leq 1 $ then 
$\beta_g(r)= c_0 \,\left[e^{(1/r)^{1/\gamma} }-1\right]^{n/2\alpha},\; r>0$.
$$
 \beta_g(r) \sim 
 \begin{cases}
 c_0 \;e^{\frac{n}{2\alpha} (1/r)^{1/\gamma} } \;\;\;{\rm as}\; \; r\rightarrow 0^+ ,\\
c_0 \; \frac{1}{r^{n/2\alpha\gamma}}   \qquad{\rm as}\;\; r\rightarrow +\infty.
\end{cases}
$$
 \item
 Let $t>0$. If $g(x)=1-e^{-tx}$
  then 
$\beta_g(r)= \frac{c_0}{t^{n/2}} \,\left[ \ln(1+\frac{1}{r-1}) \right]^{n/2},\; r>1$.
$$
 \beta_g(r) \sim 
 \begin{cases}
 \frac{c_0}{t^{n/2}} \,\left[ \ln(\frac{1}{r-1}) \right]^{n/2} 
 \;\;\,{\rm as}\; \; r\rightarrow 1^+ ,\\
\frac{c_0}{(rt)^{n/2}} 
\;\;{\rm as}\;\; r\rightarrow +\infty.
\end{cases}
$$
 \end{enumerate}
 Note that this discussion with this family of Bernstein  functions 
 is always valid when $A$ satisfies super-Poincar\'e inequality  with $\beta(t)= c_0 t^{-n/2}$.

\subsection{The Ornstein-Uhlenbeck operator}\label{ou}
 Let $A={\mathcal L}=\Delta+x.\nabla$   be the Ornstein-Uhlenbeck operator define on  $L^2(\R^n,\gamma)$
with the gaussian measure $\gamma(dx)=(2\pi)^{-n/2} e^{-\frac{\vert x\vert^2}{2}}\,dx$.
  It is well known that Gross'   logarithmic Sobolev  inequality is satisfied 
$$
\int_{\R^n} f^2\log(\vert f \vert /\vert\vert f \vert\vert_2)\,d\gamma
\leq 
( {\mathcal L}f ,f)
$$
with $( {\mathcal L}f ,f)=\int_{\R^n}
\vert\nabla f\vert^2\, d\gamma.$

 A super-Poincar\'e inequality can be deduced from Logarithmic Sobolev inequality.
  We recall the arguments.
For  any $f\in L^1(\mu)\cap L^2(\mu)$ such that $\vert\vert f \vert\vert_1=1$,  we have by Jensen's inequality:
$$
\vert\vert f \vert\vert_2^2\log \vert\vert f \vert\vert_2 
\leq
\int_{\R^n} f^2\log( \vert f \vert/\vert\vert f \vert\vert_2)\,d\gamma.
$$
By renormalization, this inequality is also satisfied   when  $\vert\vert f \vert\vert_1\leq 1$. This yields   
$$
\vert\vert f \vert\vert_2^2\log \vert\vert f \vert\vert_2 
\leq
( {\mathcal L}f ,f).
$$

Using the relation
$  xy-e^{y-1}\leq x\log x$
 for any $y\in \R$,
we deduce  
$$
\vert\vert f \vert\vert_2^2\leq t( {\mathcal L}f ,f) +  \frac{t}{2e}e^{\frac{2}{t}}\vert\vert f\vert\vert_1^2,\; t>0.
$$
On the other hand, Poincar\'e inequality deduced from Gross' inequality, trivially implies 
  $\vert\vert f \vert\vert_2^2\leq ( {\mathcal L}f ,f)+ \vert\vert f \vert\vert_1^2$. 
Together, the preceding inequalities leads to the following  formulation of super-Poincar\'e inequality
\begin{equation}\label{spg}
\vert\vert f \vert\vert_2^2\leq t( {\mathcal L}f ,f) +\beta(t) \vert\vert f\vert\vert_1^2,\; t>0,
\end{equation}
with  $\beta(t)=\frac{t}{2e}e^{\frac{2}{t}}$, $0< t\leq 1$ and  $\beta(t)=1$,  $t\geq 1$.

For instance Theorem \ref{main1} implies : 
\begin{enumerate}
\item
For any $0<\alpha<1$,
$$
\vert\vert f \vert\vert_2^2\leq t( {\mathcal L}^{\alpha}f ,f) + 
\frac{t^{\frac{1}{\alpha}}}{2e}
\,
e^{2\,t^{-1/\alpha}}
\vert\vert f\vert\vert_1^2,\quad 0<t<1,
$$
and
$$
\vert\vert f \vert\vert_2^2\leq t( {\mathcal L}^{\alpha}f ,f) + 
\vert\vert f\vert\vert_1^2,\quad t\geq 1.
$$
\item
$$
\vert\vert f \vert\vert_2^2\leq t( \log (I+{\mathcal L} )f ,f) + 
 \frac{1}{2e^3}
 \,
\frac{1}{ \left(e^{1/t}-1\right)}
 e^{2e^{1/t}}
\,
 \vert\vert f\vert\vert_1^2,\quad 0<t< \frac{1}{\log 2},
$$
and
$$
\vert\vert f \vert\vert_2^2\leq t( \log (I+{\mathcal L} )f ,f) + 
 \vert\vert f\vert\vert_1^2,\quad  t\geq  \frac{1}{\log 2}.
$$
we have 
$\beta_{\log}(t)\sim  \frac{1}{2e^3}
 \,
 e^{2e^{1/t}-1/t}
$ as $t$ goes to $0$.
\end{enumerate}
Similar inequalities can be written for the cases (iii) and (iv) considered in  Section \ref{intro}. Of course, the discussion is not limited to these cases just above.

 \section{Study of $\Psi(A)$ with $\Psi$ convex}\label{convex} 
It is  useful to deduce super-Poincar\'e or Nash inequality for $\Psi(A)$ when $A$ satisfies such inequality and
$\Psi$ is convex. The reason is that the  inverse  function of a concave  increasing function is convex and increasing.
For instance,   Bernstein functions.
In what follows, we study the following converse implication.
Assume that $g(A)$ satisfies super-Poincar\'e or Nash-type inequality then deduce a similar inequality for $A$.
\\

Let $(E_{\lambda})_{ \lambda >0}$  be the spectral resolution associated to $A$ and $\psi:[0,\infty)\rightarrow [0,\infty)$ be a measurable function.
We define $\psi(A)$ on  its domain ${\mathcal D}({\psi(A)})\subset L^2(\mu)$ as in  Section \ref{intro}.
In particular,  we have on their respective domains the following representations
$$
(Af,f)=\int_{0}^{+\infty} {\lambda}\, d(E_{\lambda}f,f),
\;\;
 (T_tf,f)=\int_{0}^{+\infty} e^{-{\lambda}t}\, d(E_{\lambda}f,f), 
\;\;
\vert\vert f \vert\vert_2^2 =\int_{0}^{+\infty} \, d(E_{\lambda}f,f).
$$ 
See \cite{scsv} Thm. 11.4.

 \begin{eprop}\label{pnashphi}
 Assume that $A$ is a   non-negative self-adjoint  
 operator satisfying Nash-type inequality (\ref{nashi}). Then for any non-negative non-decreasing 
  convex function $\Psi$ with $\Psi(0)=0$, we have
\begin{equation}\label{nashphi}
\vert\vert f\vert \vert_2^2\, (\Psi \circ D)(\vert\vert f\vert \vert_2^2)\leq (\Psi(A)f,f), 
\quad \vert\vert f\vert \vert_1\leq 1.
\end{equation}
\end{eprop}

Note that such result can be generalized  in the framework of Hilbert space  $H$ with the norm $\vert\vert f \vert\vert_1$ replaced by another control  $\Phi(f)$,  $f$ in some  subspace of  $H$,  satisfying  properties as defined in \cite{w4}.
\\

{\bf Proof}:
Write equivalently Nash-type  inequality (\ref{nashi})  as follows.
For any $f\in {\mathcal D}(A)\cap L^1(\mu)$ and  $ \vert\vert f\vert \vert_2=1$,
$$
D\left(\frac{1}
{ \vert\vert f\vert \vert_1^2 }
\right)\leq (Af,f).
$$

Since $\Psi$ is non-decreasing, we get
$$
(\Psi\circ D)\left(\frac{1}
{ \vert\vert f\vert \vert_1^2 }
\right)\leq \Psi\left[(Af,f)\right].
$$
Now, by  
functional calculus and  Jensen's inequality applied to the probability measure
$d(E_{\lambda}f,f)$, i.e. $ \vert\vert f\vert \vert_2=1$, we get 
$$
 \Psi\left[(Af,f)\right]=\Psi\left(\int_0^{+\infty} \lambda \, d(E_{\lambda}f,f)\right)
 \leq
 \int_0^{+\infty} \Psi(\lambda) \, d(E_{\lambda}f,f)= (\Psi(A)f,f).
 $$
 Thus 
 $$
(\Psi\circ D)\left(\frac{1}
{ \vert\vert f\vert \vert_1^2 }
\right)\leq 
(\Psi(A)f,f),\quad  \vert\vert f\vert \vert_2=1,\; 
f\in {\mathcal D}(A)\cap L^1(\mu).
$$
Under the assumptions on $\Psi$,  there exists $\lambda_0>0$  and  $k>0$  such that 
$k\lambda \leq \Psi(\lambda)$
 for any 
$\lambda>\lambda_0$. This implies   ${\mathcal D}(\Psi(A))\subset {\mathcal D}(A)$.
Now, by reversing the process of normalization from $L^2(\mu)$ to $L^1(\mu)$, 
it yields for any $f\in {\mathcal D}(\Psi(A))\cap L^1(\mu)$,
$$
\vert\vert f\vert \vert_2^2 \;(\Psi\circ D)\left(\frac{\vert\vert f\vert \vert_2^2}
{ \vert\vert f\vert \vert_1^2 }
\right)\leq 
(\Psi(A)f,f).
$$
Since   $\Psi\circ D$ is non-increasing,   
we deduce (\ref{nashphi}) for $ \vert\vert f\vert \vert_1\leq 1$ and  conclude the proof.
\rule{1.3mm}{2mm}
\\

In the second part of this section, we deal with the case of   super-Poincar\'e 
inequalities   generalizing  the arguments 
of  \cite{w1} used for $A^{\alpha}$, $\alpha>1$.
\begin{ethm} \label{spbspphi}
Let  $\Psi:(0,+\infty)\rightarrow (0,+\infty)$ be a non-decreasing  convex function.
Assume that
${\Psi}^*(x):=\sup_{y\in (0,+\infty)}(xy-\Psi(y))$ is a bijection from
$(0,+\infty)$  to $(0,+\infty)$
 and  that $B$ is a   non-negative  symmetric operator satisfying a super-Poincar\'e inequality with some rate function $\gamma$.
%

 Then $\Psi(B)$ satisfies also a super-Poincar\'e inequality with rate function, for $t>0$, 
$$
\gamma_{\Psi}(t)
=
\inf_{0<\varepsilon <1}\;
\frac{1}{\varepsilon} \gamma\left( \varepsilon t ({\Psi}^*)^{-1}
(\frac{1-\varepsilon}{\varepsilon t})\right).
$$
In particular,  if $\Psi(x)=x^{1/\alpha}$ with $\alpha\in (0,1)$  then $\gamma_{\Psi}(t)\leq 
\frac{1}{\alpha}
\gamma(t^{\alpha})$.
\end{ethm}

{\bf Proof}:
  Young's inequality,  for any $y$, $s>0$,
$ys\leq {\Psi}(y)+{\Psi}^*(s)$ implies for  $y=(Bf,f)$ with $f\in {\mathcal D}(B)$,  and  $s>0$,
$$
s(Bf,f)
\leq
{\Psi}((Bf,f))+
{\Psi}^*(s).
$$
Assume that $\vert\vert f\vert\vert_2=1$. By spectral representation of $B$ and by Jensen's inequality with
 $\Psi$ as convex function,
   we have  already seen that
$$
0\leq
{\Psi}((Bf,f))
\leq 
({\Psi}(B)f,f).
$$
Super-Poincar\'e implies that for any $t$, $s>0$ with $\vert\vert f\vert\vert_2=1$,
$$
1
\leq
ts(Bf,f)+\gamma(ts)\vert\vert f\vert \vert_1^2.
$$
Combining   the above inequalities, it yields for $f\in {\mathcal D}(\Psi(B))$,
$$
1
\leq
t ({\Psi}(B)f,f)+
t{\Psi}^*(s) +
\gamma(ts)\vert\vert f\vert \vert_1^2.
$$
Let $\varepsilon \in (0,1)$. Since    ${\Psi}^*$
 is a bijection,   for any fixed $t>0$,    there exists $s>0$ such that $\varepsilon= 1-t{\Psi}^*(s)$, i.e. $s=({\Psi}^*)^{-1}(\frac{1-\varepsilon}{t})$. Thus we obtain
$$
\varepsilon
\leq
t ({\Psi}(B)f,f)+
\gamma
\left(t
({\Psi}^*)^{-1}\left(
(1-\varepsilon)
 t^{-1}
 \right)
\right)
\vert\vert f\vert \vert_1^2.
$$
Changing $t$ by $\varepsilon t$ and dividing     by
$\varepsilon$,
we get 
 for any $t>0$, $\varepsilon\in (0,1)$
 and $\vert\vert f\vert \vert_2=1$,
$$
1
\leq
t ({\Psi}(B)f,f)+
\frac{1}{\varepsilon} \gamma\left( \varepsilon t ({\Psi}^*)^{-1}
\left(
\frac{1-\varepsilon}{\varepsilon t} 
\right)\right)
\vert\vert f\vert \vert_1^2.
$$
We conclude   by  changing $f$
 by $f/\vert\vert f\vert \vert_2$ and by taking the infimum over $\varepsilon\in (0,1)$.
 
 We now  prove the last statement. 
Let $\Psi(x)=x^{1/\alpha}$, we have ${\Psi}^*(s)=c_{\alpha}\,s^{\frac{1}{1-\alpha}}$
with $c_{\alpha}= (1-\alpha)\alpha^{\frac{\alpha}{1-\alpha}}$. A simple computation  yields for any $\varepsilon \in (0,1)$,
$$
\gamma_{\Psi}(t)
\leq
 \frac{1}{\varepsilon} \gamma
\left[ k_{\alpha}\,
 \varepsilon  \,
 \left(
\frac{1-\varepsilon}{\varepsilon }
 \right)^{1-\alpha}
t^{\alpha}
\right]
$$
with $k_{\alpha}=\alpha^{-\alpha} (1-\alpha)^{\alpha-1}$. Choosing $\varepsilon={\alpha}$, we conclude
$\gamma_{\Psi}(t)\leq 
\frac{1}{\alpha}
\gamma(t^{\alpha})$.
The proof is complete.
\rule{1.3mm}{2mm}
\\
 
 For the case  $g(x)=x^{\alpha}$, the function obtained in \cite{w1} is given by
${\tilde \gamma}_{\Psi}(t)=2\gamma(\frac{t^{\alpha}}{2})$. 
Since $\gamma_{\Psi}$ is usually decreasing and $\gamma(0^+)=+\infty$,   the result above  is sharper up to a multiplicative constant.
We notice  that  $\vert\vert f\vert\vert_1^2$ plays no particular role in the proof.
So, it can be replaced by some functional  
$\Phi(f)$ and $L^2(\mu)$ by a general Hilbert space as in \cite{w1}.

Now we  make the connection between Bernstein functions and convex functions. 
Assume that $g$ is a Bernstein function. Since $g$ is non-decreasing and concave,
$\Psi=g^{-1}$   is non-decreasing and  convex.
Hence, Theorem~\ref{spbspphi}  allows us to prove  a converse to
Theorem \ref{main1} about  super-Poincar\'e inequalities. 
Thus applying  Theorem~\ref{spbspphi}
 with   $B=g(A)$, we get:

\begin{ecor}\label{recip}
Let $g$ be a bijective Bernstein function and $A$ be a 
 non-negative  
 symmetric operator.  Assume that $\Psi=g^{-1}$ satisfies   the conditions of Theorem \ref{spbspphi} and  that
 that $g(A)$ satisfies a super-Poincar\'e inequality with some rate function $\gamma$.

  Then $A$  satisfies a super-Poincar\'e inequality with rate function $\gamma_{\Psi}$ given in  Theorem~\ref{spbspphi}.
%
\end{ecor}
 Corollary  \ref{recip} is  sharp in the  particular case
 $g(x)=x^{\alpha}$, $\alpha\in (0,1)$. Indeed,
assume that $A$ satisfies super-Poincar\'e inequality with rate function $\beta$. By Theorem \ref{main1}, 
$g(A)$ satisfies super-Poincar\'e   inequality  with   $\gamma=\beta_g$ given in Theorem~\ref{main1}. 
Now take $\Psi(x)=x^{1/\alpha}$ in Theorem \ref{spbspphi}, it gives back
 that $A$ satisfies super-Poincar\'e inequality with
${\tilde \beta}(t)
=
\gamma_{\Psi}(t)\leq 
\frac{1}{\alpha}
\gamma(t^{\alpha})=
\frac{1}{\alpha}\beta(t)
$.
In this case, Corollary  \ref{recip}    is  essentially an 
optimal  converse of   Theorem \ref{main1}
up to the multiplicative constant $\frac{1}{\alpha}$.

\section{Asymptotic behavior  of  $g(A)$}\label{revisit}
We briefly revisit  the relation between the asymptotic behavior of $g(A)$ and the asymptotic behavior of $A$ in terms of Poincar\'e inequality
(equivalent to a bound on the bottom of the spectrum).
The proof uses  arguments of  Theorem  \ref{main1}. We obtain Poincar\'e inequality in $L^p(\mu)$  for $g(A)$ starting from the same inequality for $A$.
See \cite{cgr} for recent results on the subject on $L^p$.
\begin{eprop}
Assume that $(T_t)_{t>0}$ is  
 a symmetric Markov semigroup   satisfying  the following inequality 
\begin{equation}\label{gap}
\vert\vert T_tf -\mu(f)\vert\vert_p
\leq
e^{-{\lambda}t} \, C(f), \quad t>0,
\end{equation}
for some $\lambda \in [0,+\infty)$, $p>1$ and some positively homogeneous functional $C(f)$.
Then  for any Bernstein function $g$ with $g(0)=0$, we have 
\begin{equation}\label{gapfrac}
\vert\vert T_t^{g} f -\mu(f)\vert\vert_p
\leq
    e^{- tg({\lambda})}\, C(f),\quad  t>0.
\end{equation}
\end{eprop}

For instance with $C(f)=\vert\vert f -\mu(f)\vert\vert_p$, see \cite{cgr}.
For $p=2$ and $C(f)=\vert\vert f -\mu(f)\vert\vert_2$, {\red \eqref{gap} } is the classical Poincar\'e inequality
$$
\vert\vert   f -\mu(f)\vert\vert_2^2\leq
\frac{1}{\lambda} 
(Af, f), \quad f\in {\mathcal D}({A}).
$$
Then similarly for  $g(A)$, we deduce the Poincar\'e inequality for $g(A)$,
$$
\vert\vert   f -\mu(f)\vert\vert_2^2\leq
\frac{1}{g(\lambda)} (g(A) f, f).
$$

{\bf Proof}:
Let $f$ such that 
$\mu(f)=0$. Then $\mu(T_sf)=\mu(T^{g}_sf)=0$, $s>0$ since  $T_s$ and $T^{g}_s$ are 
symmetric. We have for any $t>0$,
\begin{multline*}
\vert\vert T_t^{g} f -\mu(f)\vert\vert_p=
\vert\vert T_t^{g} f  \vert\vert_p=
\vert\vert
 \int_0^{+\infty} \eta_t^{g}(s)\,T_s f\,ds\vert\vert_p\\
 \leq 
  \int_0^{+\infty} \eta_t^{g}(s) \vert\vert T_s f\vert\vert_p \,ds
  \leq
       C(f) \int_0^{+\infty} \eta_t^{g}(s)   e^{-{\lambda} s}\,ds
 \leq 
    e^{- t g({\lambda})} C(f).
 \end{multline*}
This concludes the proof. \rule{1.3mm}{2mm}

\section{Functions of the Laplacian on $\R^n$}\label{fonclap}
Here, we give a direct proof of  super-Poincar\'e inequality for $g(\Delta)$ with $g$ a Bernstein function and $\Delta$ the usual Laplacian on $\R^n$.
In fact, $g$ need not be a Bernstein function. The constants are certainly not optimal.
 The proof follows the original idea used for  the Laplacian
in the paper by J. Nash \cite{N}.
\\

We use the following definition of Fourier transform
$
{\mathcal F}{f}(x)=\int_{\R^n} f(y)e^{-ixy}\, dy.
$
Let $\Delta=-\sum_{i=1}^n \frac{\partial^2}{\partial^2x_i}$. 
So,
${\mathcal F}({\Delta f})(x)
=\vert x\vert^2 {\mathcal F}{f}(x)$. The function $ {g(\Delta) f}$ is defined by its Fourier transform
$$
{\mathcal F}({g(\Delta) f})(x)=g(\vert x\vert^2){\mathcal F}{f}(x),\; x\in \R^n.
$$
The domain of $g(\Delta)$ is defined  by
$$
{\mathcal D}( g(\Delta))=\{ f\in L^2(\mu):  \int_{\R^n} \vert g(\vert x\vert^2)\vert^2\vert {\mathcal F}{f}(x)\vert^2\, dx<+\infty\}.
$$
We denote by $\omega_n $  the volume of the unit ball of $\R^n$ and set
 $c_n= (2\pi)^{-n}$.

\begin{ethm}\label{spnashgeom}
Let $g:[0,+\infty)\rightarrow[0,+\infty)$ be non-decreasing such that  $g(0)=0$.  
Let define, for $u\in[0,+\infty)$, 
$g^{\rightarrow}(u)=\sup\{s\geq 0: g(s)\leq u\}\in [0,+\infty] $.
Then we have 
\begin{enumerate}
\item
For any $t>0$ and any $f\in {\mathcal D}(g(\Delta))\cap L^1(\R^n)$, 
\begin{equation}\label{spgeom}
\vert\vert f\vert\vert_2^2\leq 
t(g(\Delta)f,f)+
{\tilde{\beta}}
\left(\frac{1}
{g^{\rightarrow}(t^{-1})}
\right)
\vert\vert f\vert\vert_1^2
\end{equation}
with ${\tilde{\beta}}(t)=c_n\, \omega_n t^{-n/2}$.
\item
For any $f\in {\mathcal D}(g(\Delta))\cap L^1(\R^n)$ with $\vert\vert f\vert\vert_1\leq 1$,  
$$
\vert\vert f\vert\vert_2^2 \, {\tilde D}_g(\vert\vert f\vert\vert_2^2)\leq  (g(\Delta)f,f)
$$
with ${\tilde  D}_g(x)=\sup_{t>0} \left(t-  \frac{t}{x} c_n\, \omega_n\left[ g^{\rightarrow}(t)\right]^{n/2}\right)$.
\end{enumerate}
\end{ethm}

When  $g$ is unbounded and $g(0)=0$, the function $g^{\rightarrow}(t)$ is well defined and finite  for any $t>0$.
The generalized inverse function  $g^{\rightarrow}$ is non-decreasing and when  $g$ is an increasing  bijection we have $g^{\rightarrow}=g^{-1}$. 
The function ${\tilde{\beta}}_g(t):={\tilde{\beta}}
\left(\frac{1}
{g^{\rightarrow}(t^{-1})}
\right)$ in (\ref{spgeom}) is similar to $\beta_g(t)$ in Theorem~\ref{main1} when $g$ is invertible.
If we assume that $g$  is bounded then   $g(\Delta)$ is a bounded operator
and  ${\tilde{\beta}}_g(t)=+\infty$  when $t \leq 1/\vert\vert g\vert\vert_{\infty}$.
In that case,  the inequality  (\ref{spgeom}) is meaningful only for $t>1/\vert\vert g\vert\vert_{\infty}$. 
Note that this restriction already appears in Theorem~\ref{main1} when $g$ is bounded.
\\

{\bf Proof}:  
Let   $f \in {\mathcal D}( g(\Delta))\cap L^1(\mu)$ and $t>0$.
By Plancherel formula,
\begin{multline*}
\vert\vert f\vert\vert_2^2=c_n\int_{\R^n} \vert {\mathcal F} f\vert ^2(x)\, dx
=
c_n\int_{\{ x\in \R^n: 1\,\leq \, tg(\vert x\vert^2)\}}\!\!\!\!\!\!\!\!\!\!\!\!\! \vert {\mathcal F} f\vert ^2(x)\, dx
+
c_n\int_{\{ x\in \R^n: 1\,>\,tg(\vert x\vert^2)\}}\!\!\!\!\!\!\!\!\!\!\!\!\! \vert {\mathcal F} f\vert ^2(x)\, dx\\
\leq 
c_n \,t \int_{\R^n} 
 g(\vert x\vert^2) \vert {\mathcal F} f\vert ^2(x)\, dx
+
c_n\, 
\vert\vert  {\mathcal F} f\vert\vert_{\infty}^2 \,
V(\{ x\in \R^n: g(\vert x\vert^2)< \frac{1}{t}\})
\end{multline*}
where $V(\Omega)$ is the Lebesgue measure of the set $\Omega\subset \R^n$.
Now, since $g(r)\leq  u$ implies $r\leq g^{\rightarrow}(u)$
and
$\vert\vert  {\mathcal F} f\vert\vert_{\infty}\leq  \vert\vert f\vert\vert_1$,
we deduce  for any $t>0$,
 $$
 \vert\vert f\vert\vert_2^2 
 \leq   t\,(g(\Delta f ,f)
 +
c_n\, \vert\vert  f
\vert\vert_1^2  \,V\left(\{x\in \R^n: \vert x\vert\leq \sqrt{g^{\rightarrow}(t^{-1})}\}\right).
$$
 This concludes (i).

We prove the second part by applying   Proposition \ref{nashpoin}
with $\tilde\beta(t)= \, c_n \, \omega_n   
 \left[g^{\rightarrow}(t^{-1})\right]^{n/2}
$. This completes the proof.
 \rule{1.3mm}{2mm}

\section{Appendix on Legendre transform}\label{appen2}

In the first part of this section,  we are interested to discuss  the properties of $D$ directly from the properties of $\beta$  independently  of the set ${\mathcal G}$ defined in  Proposition \ref{nashpoin}. In the second part of this section,  reversing the role of $\beta$
and $D$ leads to a similar discussion. The conditions introduced  here are usually satisfied in the applications. 

\begin{elem}\label{betad}
Let $\beta $ be a non-negative function  on $(0,+\infty)$ and set
$$
D(x)=\sup_{t>0} \left\{t-\frac{t}{x}\beta(1/t)\right\}\in (-\infty,+\infty],\quad x>0.
$$ 
\begin{enumerate}
\item
If $\lim_{t\rightarrow 0^+} t\beta(1/t)=0$ then $D$ is non-negative. 
This condition is satisfied if $\beta$ is bounded above at infinity, 
in particular if $\beta$ is non-increasing.
\item
If  $\beta(0^+)=+\infty$ then $D(x)$ is finite for any $x>0$. Moreover, the function  $x\rightarrow xD(x)$ 
is convex, non-decreasing on $(0,+\infty)$ and $D$ is  continuous.
\end{enumerate}
\end{elem}
{\bf Proof}:
\begin{enumerate}
\item
 Assume that $\lim_{t\rightarrow 0^+} t\beta(1/t)=0$ then $D(x)\geq \lim_{t\rightarrow 0^+} (t-\frac{t}{x}\beta(1/t))=0$
for any $x>0$.
Obviously, if $\beta$ is   bounded above at infinity (e.g.  non-increasing) and  non-negative then  $\lim_{t\rightarrow 0^+} t\beta(1/t)=0$.
\item
Assume $\beta(0^+)=+\infty$.
Fix $x>0$ then  there exists   $t_x>0$   such that for any $t>t_x$, $x< \beta(1/t)$. 
So,  $t-\frac{t }{x}\beta(1/t)<0$ when $t>t_x$. 
If $0< t\leq t_x$
then $t-\frac{t }{x}\beta(1/t)\leq t_x$  since $\beta$ is non-negative. Therefore  $D(x)\leq t_x$  and  $D(x)$  is finite. 
Now, the function  $ x\rightarrow h^*(x)=xD(x)=\sup_{t>0}\left( tx-t\beta(1/t)\right)$ is     convex on $(0,+\infty)$.  Consequently,  $h^*$ and $D$ are continuous.
\end{enumerate}
 The proof is complete.
 \rule{1.3mm}{2mm}

Now, we study properties of $\beta$ in terms of $D$  defined as in Theorem \ref{main2}.
 Natural conditions on $D$ comes from the previous lemma. The discussion is similar.
\begin{elem}\label{dbeta}
Let $D:(0,+\infty)\rightarrow \R$ be a fixed function and set
$
 \beta(r)=\sup_{x>0}\!\{ x-rxD(x)\}.
 $
\begin{enumerate}
\item
If $\lim_{x\rightarrow +\infty} D(x)=+\infty$ and $D$ is non-negative  then $\beta(r)$ is finite for  any $r>0$,   convex, continuous and non-increasing.
\item
If $\lim_{x\rightarrow 0} xD(x)=0$ 
then $\beta$ is non-negative.
\end{enumerate}
\end{elem}
The proof is   similar to  Lemma~\ref{betad}. 
Indeed, let $\beta$ and $D$ as in Lemma~\ref{betad}.
Note that $h^*(x):=xD(x)=\sup_{t>0} (tx-t\beta(1/t))$ is the Lengendre transform 
(or complementary function) of $h(t)=t\beta(1/t)$ thus the afferent theory applies.
See  for instance  \cite{rr} p.6 for a discussion about Young functions and  p.13
 for the specific class  of N-functions.
 Usually   $h$ is obtained  from $h^*$ by the same   formula, i.e.
$h(t)= \sup_{x>0} (tx-h^*(x))$, $t>0$. 
In that case, we recover the definition of $\beta$ in terms of $D$ in 
Lemma  \ref{dbeta} by the formulas  $\beta(t)=t h(1/t)$ and $D(x)= h^*(x)/x$.
Here are some examples of couples of  N-functions. 
 They appear 
as asymptotics of  functions $\beta$ or $D$ of our  examples  in   Section \ref{exset}. 
Let   $1<p,q<+\infty$ with $1/p+1/q=1$. 

\begin{enumerate}
\item
$( h_1(t), h_1^*(x))=(t^p/p, x^q/q)$.
\item
$(h_2(t), h_2^*(x))=(e^t-t-1,(1+x)\ln(1+x)-x)$.
\item
$(h_3(t), h_3^*(x))=((1+t)\ln(1+t)-t,e^x -x-1)$.
\item
 $h_4(t)=e^{t^p}-1$, $h_4^*$:  no explicit form.
\end{enumerate}
But one can prove that
$h_4^*(x)\sim x\left( \ln x \right)^{1/p}$ as $x\rightarrow +\infty$
and
$h_4^*(x)\sim c_q \,  x^q $ as $x\rightarrow 0^+$ with $1/p+1/q=1$ and $c_q=(p-1)\left(\frac{1}{p}\right)^q$.

Of course, in applications, functions like $c_1h(c_2 t)$, $c_i>0$  should be considered or functions having asymptotics of this type.
Indeed, in practice $h$ is not  exactly an N-function  but often close to  such function. Fortunately,  it doesn't cause much trouble in practice.
This justifies  the interest of both propositions just above.
Now,  we mentioned some cases  where the asymptotics of these functions really  appear in our applications.
Recall the relations  $\beta(t)=th(1/t)$ and $D(x)= h^*(x)/x$.

\begin{enumerate}
\item
For $h_1(t)=c_0\,t^p$ then we get $\beta_1(t)=\frac{c_1}{t^{\nu}}$ where $\nu=p-1>0$
and $D(x)=c_{2}\, x^q$  with $1/p+1/q=1$ and $q=1+\frac{1}{\nu}$.
Such cases correspond in Section \ref{eucli}  to  the fractional Laplacian $\Delta^{\alpha}$ on the Euclidean space and   to 
Section \ref{hypo}   with $L^{\alpha}$ where  $L$ is a sum of vector fields satisfying H\"ormander's condition and  $\nu=\frac{n}{2\alpha}$.
For all these cases: $D(x)=c_3\, x^{\frac{2\alpha}{n}}$.
\item
The function  $h_2(t)=e^t-t-1$ leads to $\beta_2(t)\sim te^{\frac{1}{t}}$ as $t\rightarrow 0^+$
and $\beta_2(t)\sim \frac{1}{2t}$ as $t\rightarrow +\infty$. This situation  is realized 
up to multiplicative constants 
by the Ornstein-Uhlenbeck operator in Section \ref{ou} as far as the local behavior, i.e. $t\rightarrow 0^+$,
is concerned.
In that case, $D(x)\sim  \ln x$ as $x\rightarrow +\infty$.
\item
The function $h_3(t)=(1+t)\ln(1+t)-t$ gives $\beta_3(t)\sim \ln (1/t)$
as  $t\rightarrow 0^+$.
 Using  results of Section \ref{fonclap}   with $n=2$,  we set $A=g(\Delta)\geq  0$ 
 with $g(y)=e^{y/4\pi}-1,\; y\geq 0$. This
 provides an example of positive operator such that the super-Poincar\'e inequality (\ref{spgeom}) is satisfied with $\beta(t)= \ln (1+ \frac{1}{4\pi^2 t}),\; t>0$ and $D(x)\sim \frac{4\pi^2}{x}e^{x-1}$ as $x\rightarrow +\infty$. We also have
 $D(x)\sim \pi^2x$ as $x\rightarrow 0^+$.
\item
For the function   $h_4(t)=e^{t^p}-1$ with $1<p<+\infty$,  we deduce $\beta_4(t)\sim t e^{\frac{1}{t^p}}=:\tilde {\beta}(t)$ as $t\rightarrow 0^+$ and $e^{(1-\frac{1}{p}) \frac{1}{t^p}}\leq  \tilde {\beta}(t) 
\leq e^{\frac{1}{t^p}}$  when  $t\in (0,1)$.
Examples with such behavior are  given in the   Riemannian setting of Section \ref{riem}  with 
$p=\frac{\delta}{2(\delta-1)}$
where  $1<\delta <2$ in  (\ref{exporiem}).
In that  case,  $D(x)\sim \left(\ln x\right)^{1/p}$ as $x\rightarrow +\infty$. 
Note that it is  a general fact that 
 the behavior of $\beta(t)$ as $t\rightarrow 0^+$ determines the behavior of $D(x)$ as $x\rightarrow +\infty$
 and conversely.
\end{enumerate}

It will be interesting to know if there exists an operator $A$ satisfying super-Poincar\'e inequality 
with $\beta(t)\sim t h_3(1/t)\sim \ln (1/t)$ as  $t\rightarrow 0^+$ and  $(Af,f)$ a Dirichlet form. 
\\

{\em Acknowlegments:} 
The authors    thank the anonymous referee  for   helpful comments
and suggestions to improve the paper.
This research was supported in part by the ANR project EVOL. 
The second author thanks the CNRS for a period of delegation  
during which this paper
has been completed. 
\\

{\em Note}: The article  \cite{w1}   can be obtained on demand to the author and difficult by other sources.

\section*{References}

\bibitem[B-M]{bm}Bendikov A.D. and  Maheux P.:
{\em  Nash type inequalities for fractional powers of non-negative self-adjoint operators }. 
 Trans. Amer. Math. Soc.  359, no. 7  (2007), 3085-3097.

\bibitem[B-F]{bf}Berg Ch. and  Forst G.: 
{\em Potential theory on locally compact
Abelian groups.} Erg. der Math. und ihrer Grenzgeb., Band 87, Springer-Verlag
1975.

\bibitem[Bi-M]{bim}Biroli M. and Maheux P.:
{\em Super Logarithmic Sobolev inequalities and Nash-type inequalities for sub-markovian 
symmetric semigroups.}{\verb"[hal-00465177, v1]"}.

\bibitem[C-L]{cl}Carlen E.A. and Loss M.:
{\em Sharp constant in Nash's inequality}. 
Internat. Math. Res. Notices 1993, no. 7, 213--215. 

\bibitem[C-G-R]{cgr}Cattiaux  P.,  Guillin A. and   Roberto C.:
 {\em Poincar\'e inequality and the $L^p$  convergence of semi-groups.} 
   Electronic Communications in Probability, 2010, Vol.15,	
 270--280.
 
\bibitem[C]{c}Coulhon T.:
{\em Ultracontractivity and Nash type inequalities}
J.Funct.Anal.141 (1996), p.510-539.

\bibitem[D]{d}Davies  E.B.:
 {\em Heat kernels and spectral theory}. 
 Cambridge Tracts in Mathematics, 92. Cambridge University 
Press, Cambridge, 1989. 

\bibitem[G-H]{gh} 
Grigor'yan,A.,
Hu, J.: {\em Upper bounds of heat kernels on doubling spaces.}
Preprint, 2010, 
available at \verb" http://www.math.uni-bielefeld.de/~grigor/pubs.htm".

\bibitem[J]{j}Jacob N.:
{\em Pseudo-Differential Operators and Markov Processes}. Vol. 1: Fourier Analysis and Semigroups. Imperial College Press, London, 2001.

\bibitem[M1]{m1}Maheux P.:
{\em New Proofs of Davies-Simon's Theorems about Ultracontractivity and Logarithmic Sobolev Inequalities related to Nash Type Inequalities.} Preprint available on
 {\verb"ArXiv: math/0609124".}

\bibitem[N]{N} Nash J.:
{\em Continuity of solutions of parabolic and elliptic equations}.
Amer.J. Math. 80, 1958,pp.931-954.

\bibitem[R-R]{rr}Rao M. M. and  Ren Z. D.: 
{\em Theory of Orlicz spaces}, Monographs and Textbooks in Pure 
and Applied Mathematics, vol. 146, Marcel Dekker Inc., New York, 1991. 
 
  \bibitem[Sc-S-V]{scsv}
 Schilling R.L., Song R. and Vondra\v{c}ek Z.:
 {\em  Bernstein functions. Theory and applications}.
  de Gruyter Studies in Mathematics, 37. 
  Walter de Gruyter and Co., Berlin, 2010.
 
 \bibitem[S-S-V]{ssv}\v{S}iki\'c  H., Song R. and  Vondra\v{c}ek Z.:
  {\em 
  Potential theory of geometric stable processes}.
    Probab. Theory Related Fields  135  (2006),  no. 4, 547-575. 
 
\bibitem[VSC]{vsc}Varopoulos N.T., Saloff-Coste  L. and Coulhon, T.: {\em Analysis
and Geometry on Groups.} Cambridge University Press (1992).

\bibitem[W1]{w1}Wang F.-Y.: 
{\em Functional Inequalities for Dirichlet Operators with Powers}.
Chinese Sci.
Tech. Online, 2007.

\bibitem[W2]{w2}Wang F.-Y.:
 {\em  Functional inequalities for empty essential spectrum}.
  J. Funct. Anal.  170  (2000),  no. 1, 219--245. 

\bibitem[W3]{w3}Wang, F.-Y.: 
 {\em Functional inequalities for the decay of sub-Markov semigroups}. 
Potential Analysis, 18 (2003), no. 1, 1--23. 

 \bibitem[W4] {w4}Wang, F.-Y.:
{\em Functional Inequalities, Markov Processes and Spectral Theory}.
 Science Press, Beijing, 2004.

\end{document}